\documentclass[reqno]{amsart}  
\pdfoutput=1
\usepackage{hyperref}
\usepackage{graphicx}
\usepackage{enumerate} 
\usepackage{upref}
\usepackage{verbatim} 
\usepackage{color}

\newcommand*{\mailto}[1]{\href{mailto:#1}{\nolinkurl{#1}}}

\newtheorem{theorem}{Theorem}[section]
\newtheorem{lemma}[theorem]{Lemma}

\numberwithin{equation}{section}
\newtheorem{definition}[theorem]{Definition}
\allowdisplaybreaks

\newcommand{\act}{\bullet}

\newcommand{\Wperun}{{W^{1,1}_{\text{\rm per}}}}
\newcommand{\Wper}{{W^{1,\infty}_{\text{\rm per}}}}

\renewcommand{\H}{\mathcal{H}}

\newcommand{\Gr}{G}

\newcommand{\sj}{\sum_{j=1}^{n}}
\newcommand{\PP}{\mathcal{P}}
\newcommand{\Q}{\mathcal{Q}}

\newcommand{\dott}{\,\cdot\,}

\newcommand{\D}{\ensuremath{\mathcal{D}}}

\newcommand{\F}{\ensuremath{\mathcal{F}}}

\newcommand{\Real}{\mathbb{R}}

\newcommand{\muacZ}{\mu_{0, \text{\rm ac}}}

\DeclareMathOperator{\id}{id}

\DeclareMathOperator{\sgn}{sgn}

\newcommand{\beq}{\begin{equation}}
  \newcommand{\eeq}{\end{equation}}
\newcommand{\bal}{\begin{align}}
  \newcommand{\eal}{\end{align}}

\newcommand{\norm}[1]{\left\Vert#1\right\Vert}
\newcommand{\abs}[1]{\left\vert#1\right\vert}

\numberwithin{equation}{section}

\begin{document}

\title[Periodic 2CH system]{Periodic conservative solutions for the two-component Camassa--Holm system}

\author[K. Grunert]{Katrin Grunert}
\address[K. Grunert]{\newline 
  Department of Mathematical Sciences\\ Norwegian University of Science and Technology\\ NO-7491 Trondheim\\ Norway}
\email{\mailto{katring@math.ntnu.no}}
\urladdr{\url{http://www.math.ntnu.no/~katring/}}

\author[H. Holden]{Helge Holden}
\address[H. Holden]{\newline
  Department of Mathematical Sciences\\
  Norwegian University of Science and Technology\\
  NO-7491 Trondheim\\ Norway\newline {\rm and} \newline Centre of
  Mathematics for Applications\\ University of Oslo\\
  NO-0316 Oslo\\ Norway}
\email{\mailto{holden@math.ntnu.no}}
\urladdr{\url{http://www.math.ntnu.no/~holden/}}

\author[X. Raynaud]{Xavier Raynaud}
\address[X. Raynaud]{\newline
  Centre of Mathematics for Applications\\
  University of Oslo\\ NO-0316 Oslo\\ Norway}
\email{\mailto{xavierra@cma.uio.no}}
\urladdr{\url{http://folk.uio.no/xavierra/}}

\date{\today} 
\thanks{(Research supported in part by the Research Council of Norway, and the Austrian Science Fund (FWF) under Grant No.~J3147.}  
\thanks{In: {\em Spectral Analysis, Differential Equations and Mathematical Physics}, Proc. Symp. Pure Math., Amer. Math. Soc. (to appear)}

\dedicatory{Dedicated with admiration to Fritz Gesztesy on the occasion of his sixtieth anniversary}

\begin{abstract}
  We construct a global continuous semigroup of weak  periodic conservative solutions to the two-component Camassa--Holm system, $u_t-u_{txx}+\kappa u_x+3uu_x-2u_xu_{xx}-uu_{xxx}+\eta\rho\rho_x=0$ and  $\rho_t+(u\rho)_x=0$, for initial data $(u,\rho)|_{t=0}$ in $H^1_{\rm per}\times L^2_{\rm per}$. It is necessary to augment the system with an associated energy to identify  the conservative solution. We study the stability of these periodic solutions by constructing a Lipschitz metric. Moreover, it is proved that if the density $\rho$ is bounded away from zero, the solution is smooth. Furthermore, it is shown that  given a sequence $\rho_0^n$ of initial values for the densities that tend to zero, then the associated solutions $u^n$ will approach the global conservative weak solution of the Camassa--Holm equation. Finally it is established how the characteristics govern the smoothness of the solution.
\end{abstract}
\maketitle

\section{Introduction}
In this paper we analyze periodic and conservative weak global solutions of the two-component 
Camassa--Holm (2CH) system   which reads (with $\kappa\in\Real$ and $\eta\in(0,\infty)$)
\begin{subequations}\label{eq:chkappa00}
  \begin{align}
    u_t-u_{txx}+\kappa u_x+3uu_x-2u_xu_{xx}-uu_{xxx}+\eta\rho\rho_x&=0,\\
    \rho_t+(u\rho)_x&=0.
  \end{align}
\end{subequations}
The special case when $\rho$ vanishes identically reduces the system to the celebrated and well-studied Camassa--Holm (CH) equation, first studied in the seminal paper \cite{CH:93}. The present system was first introduced by Olver and Rosenau in \cite[Eq.~(43)]{OlverRosenau}, and derived in the context of water waves in \cite{MR2474608}, showing  $\eta$ positive and $\rho$ nonnegative to be the physically relevant case. Conservative solutions  on the full line for the 2CH system have been studied, see, e.g., \cite{GHR2}. However, periodic and conservative solutions for the 2CH system have not been analyzed so far, and this paper aims to fill that gap. It offers some technical challenges that will be described below.

The 2CH system can suitably be rewritten as  
\begin{subequations}
  \label{eq:rewchsys10}
  \begin{align}
    \label{eq:rewchsys11}
    u_t+uu_x+P_x&=0,\\
    \label{eq:rewchsys12}
    \rho_t+(u\rho)_x&=0,
  \end{align}
\end{subequations}
where $P$ is implicitly defined by
\begin{equation}
  \label{eq:rewchsys13}
  P-P_{xx}=u^2+\kappa u+\frac12u_x^2+\eta\frac12\rho^2.
\end{equation}

The reason for the intense study of the CH equation is its  surprisingly rich structure. In the context of the present paper, the focus  is on the wellposedness of global weak solutions of the Cauchy problem. There is an intrinsic dichotomy in the solution that appears after wave breaking, namely between solutions characterized  either by conservation or dissipation of the associated energy. The two classes of solutions are for obvious reasons denoted conservative and dissipative, respectively.  The fundamental nature of the problem can be understood by the following pregnant example, for simplicity presented here on the full line, rather than the periodic case. The CH equation with $\kappa=0$ has as special solutions so-called multipeakons given by
\begin{equation*}
  u(t,x)=\sum_{i=1}^n p_i(t)e^{-\abs{x-q_i(t)}},
\end{equation*}
where the $(p_i(t), q_i(t))$ satisfy the explicit system of ordinary differential equations
\begin{equation*}
  \dot q_i=\sj p_je^{-\abs{q_i-q_j}},\quad
  \dot p_i=\sj p_ip_j\sgn(q_i-q_j)e^{-\abs{q_i-q_j}}.
\end{equation*}
In the special case of $n=2$ and $p_1=-p_2$ and $q_1=-q_2<0$ at $t=0$, the solution consists of two ``peaks'', denoted peakons,  that approach each other. At time $t=0$  the two peakons annihilate each other, an example of  wave breaking, and the solution satisfies $u=u_x=0$ 
pointwise at that time.  For positive time two possibilities exist; one is to let the solution remain equal to zero (the dissipative solution), and other one being that that two peakons reemerge (the conservative solution).  A more careful analysis reveals that the  $H^1(\Real)$ norm of $u$ remains finite, while $u_x$ becomes singular, at $t=0$, and there is an accumulation of energy in the form of a Dirac delta-function at the point of annihilation.  The consequences  for the wellposedness of the Cauchy problem are severe. The continuation of the solution past wave breaking has been studied, see \cite{BreCons:07, BreCons:09,HolRay:07,HolRay:09}.  The method to handle the dichotomy is by reformulating the equation in Lagrangian variables, and analyze carefully the behavior in those variables.  We will detail this construction later in the introduction. 

The 2CH system has, in spite of its brief history, been studied extensively, and it is not possible to include a complete list of references here. However, we mention \cite{WangHuangChen,GHR2}, where a similar approach to the present one, has been employed. The case with $\eta=-1$ has been discussed in 
\cite{eschlechyin:07}; our approach does not extend to the case of $\eta$ negative.  In \cite{GuanYin2010a} it is shown that if the initial density $\rho_0>0$, then the solution exists globally and this result is extended here to a local result, Theorem \ref{th:presreg}, where we show how the characteristics govern the local smoothness. 
For other related results pertaining to the present system, please see \cite{GuanYin2010a,GuiLiu2011,GuiLiu2010}.
There exists other two-component generalizations of the CH equation than the one studied here; see, e.g., \cite{ChenLiu2010, FuQu:09,GuanYin2010,GuoZhou2010,Kuzmin}.

We now turn to the discussion of the present paper. For simplicity we assume  that
$\eta=1$ and 
$\kappa=0$.  We first make a change from Eulerian to Lagrangian variables and
introduce a new energy variable. The change of variables, which we now will detail,  is related to the one used in \cite{HolRay:07} and, in particular, \cite{GHR:12}.  Assume that $(u,\rho)=(u(x,t),\rho(x,t))$ is a solution of \eqref{eq:chkappa00}, and define the characteristics $y=y(t,\xi)$ by
\begin{equation*}
  y_t(t,\xi)=u(t,y(t,\xi))
\end{equation*}
and the Lagrangian velocity by
\begin{equation*}
  U(t,\xi)=u(t,y(t,\xi)).
\end{equation*} 
By introducing the Lagrangian energy  density $\nu$ and density $r$ by
\begin{equation*}
  \begin{aligned}
    \nu(t,\xi)&=u^2(t,y(t,\xi))y_\xi(t,\xi)+u_x^2(t,y(t,\xi))y_\xi(t,\xi)+\rho^2(t,y(t,\xi))y_\xi(t,\xi), \\
    r(t,\xi)&=\rho(t,y(t,\xi))y_\xi(t,\xi),
  \end{aligned}
\end{equation*}
we find that the system can be rewritten as
(introducing $\zeta(t,\xi)=y(t,\xi)-\xi$ for technical reasons)  
\begin{subequations}\label{sys:persys0}
  \begin{align}
    \zeta_t&=U,\\
    U_t&=-Q,\\
    \nu_t&=-2QUy_\xi+(3U^2-2P)U_\xi,\\
    r_t&=0,
  \end{align}
\end{subequations}
where the functions $P$ and $Q$ are explicitly given by \eqref{eq:Psimp1} and  \eqref{eq:Qsimp1}, respectively.  We then establish the existence of a unique global solution for this system (see Theorem~\ref{th:global}), and  we show that the solutions form a continuous semigroup in an appropriate norm.  
In order to solve the Cauchy problem \eqref{eq:rewchsys10} we have to choose the initial data appropriately. To accommodate for the possible concentration of energy we augment the natural initial data $u_0$ and $\rho_0$ with a nonnegative Radon measure $\mu_0$ such that the absolutely continuous part $\muacZ$   equals 
$\muacZ=(u_0^2 +u_{0,x}^2+\rho_0^2)\,dx$.  The precise translation of these initial data is given in Theorem \ref{defL}.  One then solves the system in Lagrangian coordinates. The translation back to Eulerian variables is described in Definition~\ref{def:FtoE}. However, there is an intrinsic problem in this latter translation if one wants a continuous semigroup. This is due to the problem of \textit{relabeling}; to each solution in Eulerian variables there exist several distinct solutions in Lagrangian variables as there are additional degrees of freedom in the Lagrangian variables. In order to resolve this issue to get a continuous semigroup, one has to identify Lagrangian  functions corresponding to one and the same Eulerian solution. This is treated in Theorem~\ref{thm:LMPI}.  The main existence theorem, Theorem~\ref{th:main2},  states that for  $u_0\in H^1_{\rm per}$ and $\rho_0\in L^2_{\rm per}$ and $\mu_0$ a nonnegative Radon measure with absolutely continuous part $\muacZ$   such that $\muacZ=(u_0^2+u_{0,x}^2+\rho_0^2)\,dx$, there exists a continuous semigroup $T_t$ such that $(u,\rho)$, where $(u,\rho,\mu)(t)=T_t(u_0,\rho_0,\mu_0)$, is a weak global and conservative solution of the 2CH system. In addition, the measure $\mu$ satisfies
\begin{equation*}
  (\mu)_t+(u\mu)_x=(u^3-2Pu)_x,
\end{equation*}
weakly. Furthermore, for almost all times the measure $\mu$ is absolutely continuous and 
$\mu=(u^2+u_{x}^2+\rho^2)\,dx$.  

The solution so constructed is not Lipschitz continuous in any of the natural norms, say $H^1$ or $L^p$. Thus it is an intricate problem to identify a metric that deems the solution Lipschitz continuous, see \cite{GHR,GHRb:10}. For a discussion of Lipschitz metrics in the setting of the Hunter--Saxton equation and relevant examples from ordinary differential equations, see \cite{BHR}. The metric we construct here has to distinguish between conservative and dissipative solutions, and it is closely connected with the construction of the semigroup in Lagrangian variables.  We commence by defining a metric in Lagrangian coordinates. To that end, let
\begin{equation}
  \label{eq:defJ0}
  J(X_\alpha,X_\beta)=\inf_{f,g\in\Gr}\norm{X_\alpha\act f-X_\beta\act g}_E.
\end{equation}
Here $G$ contains the  labels used for the relabeling, see Definition \ref{def:group}, and $X\act f$ denotes the solution $X$ with label $f$.
The function $J$ is invariant with respect to relabeling, yet it is not a metric as it does not satisfy the triangle inequality.   Introduce
$d(X_\alpha,X_\beta)$  by 
\begin{equation}
  \label{eq:defdist0}
  d(X_\alpha,X_\beta)=\inf \sum_{i=1}^NJ(X_{n-1},X_n), \quad X_\alpha,X_\beta\in\F,
\end{equation}
where the infimum is taken over all finite sequences
$\{X_n\}_{n=0}^N\in\F$ satisfying
$X_0=X_\alpha$ and $X_N=X_\beta$.  This will be proved to be a Lipschitz metric in Lagrangian variables. Next the metric is transformed into Eulerian variables, and 
Theorem~\ref{th:main} identifies a metric, denoted $ d_{\D^M}$, such that the solution is Lipschitz continuous.

Due to the non-local nature of $P$ in  \eqref{eq:rewchsys13}, see \eqref{eq:nonlocal}, information travels with infinite speed. Yet, we show in Theorem \ref{th:presreg}  that regularity is a local property in the following precise sense.   A solution is said to be $p$-regular, with $p\geq 1$ if 
\begin{equation*}
  u_0\in W^{p,\infty}(x_0,x_1), \quad \rho_0\in W^{p-1,\infty}(x_0,x_1), \quad \text{and } \mu_0=\mu_{0,ac} \text{ on } (x_0,x_1), 
\end{equation*}
and that 
\begin{equation*}
  \rho_0(x)^2\geq c>0
\end{equation*}
for $x\in(x_0,x_1)$. If the initial data $(u_0,\rho_0,\mu_0)$ is $p$-regular, then the solution
$(u,\rho,\mu)(t,\dott)$, for $t\in \Real_+=[0,\infty)$,  remains $p$-regular on the interval $(y(t,\xi_0),y(t,\xi_1))$, where $\xi_0$ and $\xi_1$ satisfy $y(0,\xi_0)=x_0$ and $y(0,\xi_1)=x_1$ and are defined as 
\begin{equation*}
  \text{$\xi_0=\sup\{\xi\in\Real\ |\ y(0,\xi)\leq x_0\}$  and 
    $\xi_1=\inf\{\xi\in\Real\ |\ y(0,\xi)\geq x_1\}$.}
\end{equation*}

\begin{figure}
  \centering
  \includegraphics[width=12cm]{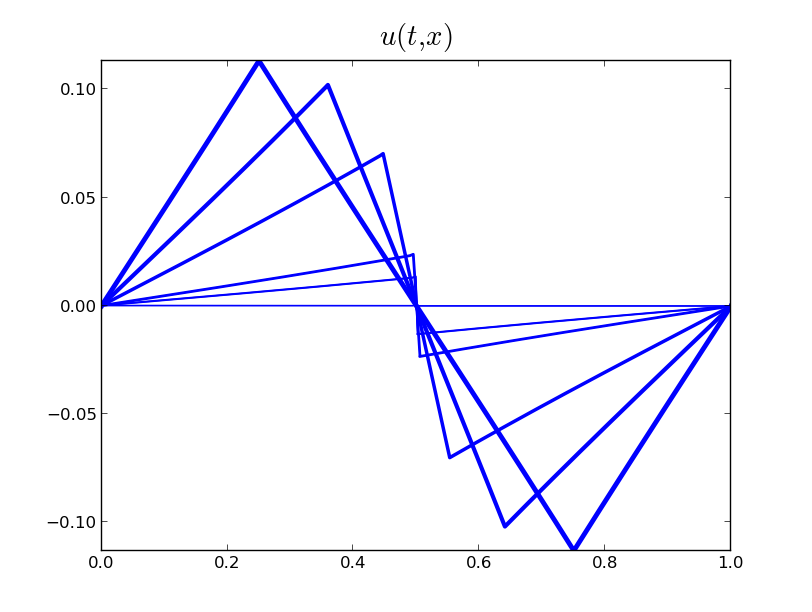} \\
  \includegraphics[width=12cm]{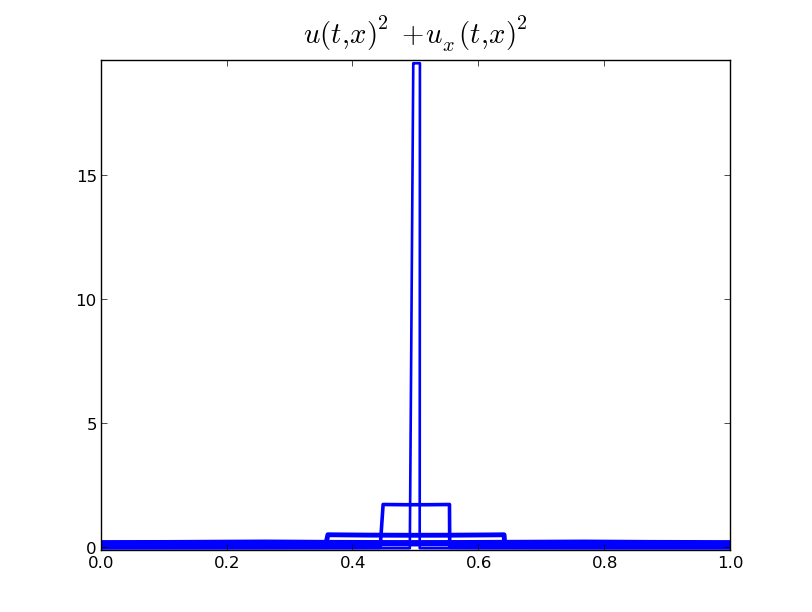}
  \caption{Plot of a peakon anti-peakon solution $u$ of the CH equation
    with $\rho$ identically zero at all times (the
    thinner the curve is, the later time it
    represents). In this case, we obtain a
    conservative solution of the scalar
    Camassa--Holm equation. We observe that the
    total energy $u^2+u_x^2$ converges to a multiple of a 
    Dirac delta function at $t=t_c$.}
  \label{fig:coll}
\end{figure}

\begin{figure}
  \centering
  \includegraphics[width=9cm]{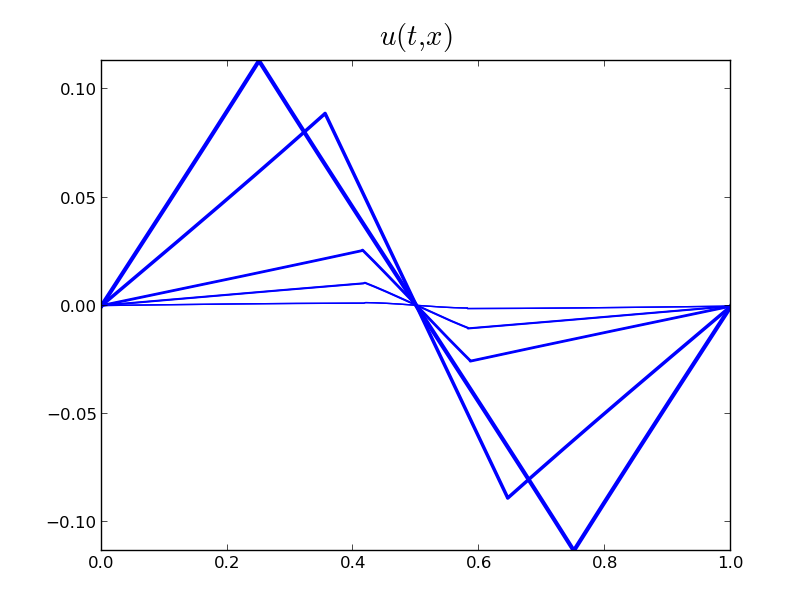}\\ 
  \includegraphics[width=9cm]{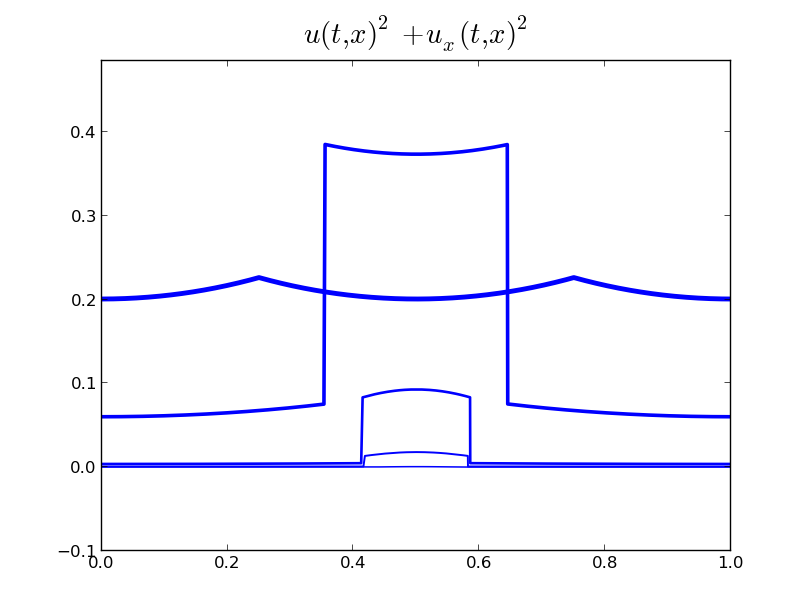}\\  \includegraphics[width=9cm]{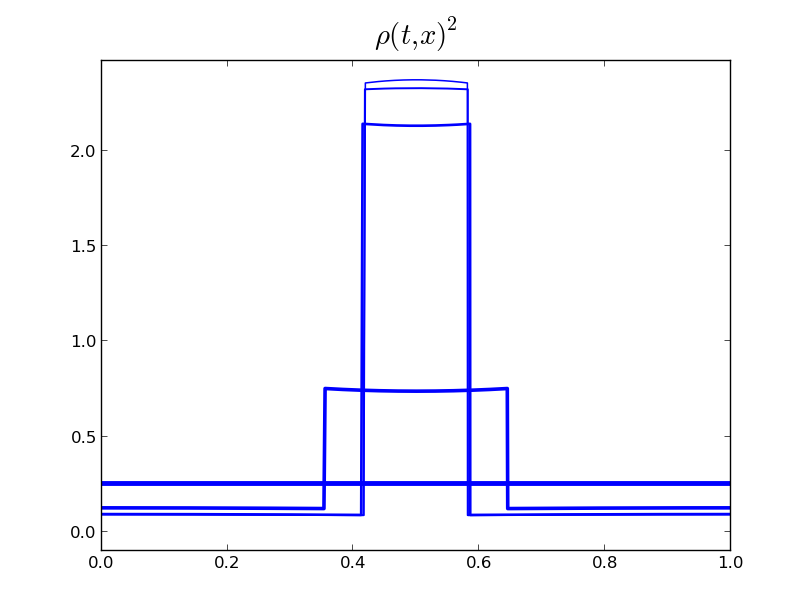} 
  \caption{Here we employ the same initial condition as in Figure~\ref{fig:coll} for $u$ while
    $\rho_0(x)=0.5$.  The total energy
    $\int_0^1(u^2+u_x^2+\rho^2)\,dx$ is
    preserved. We observe first a concentration of
    the part of the energy given by $u^2+u_x^2$. However, as we get closer to $t_c$, there is a transfer of energy from $u^2+u_x^2$ to $\rho^2$.}
  \label{fig:col2}
\end{figure}

It is interesting to consider how the standard
CH equation is obtained when the density $\rho$
vanishes since the CH equation formally is
obtained when $\rho$ is identically zero in
the 2CH system. In order to analyze the behavior of the
solution, we need to have a sufficiently strong
stability result.  Consider a sequence of
initial data $(u_0^n,\rho_0^n, \mu_0^n)$ such
that $u_0^n\to u_0$ in $H^1_{per}$,
$\rho_0^n\to 0$ in $L^2_{per}$ with
$\rho_0^n\ge d_n>0$ for all $n$. Assume that the
initial measure is absolutely continuous, that
is,
$\mu_0^n=\muacZ^n=((u_{0,x}^n)^2+(\rho_0^n)^2)\,dx$. Then
we show in Theorem \ref{th:approxCH} that the
sequence $u^n(t)$ converges in $L^\infty_{per}$
to the weak, conservative global solution of the
Camassa--Holm equation with initial data
$u_0$. To illustrate this result we have plotted, in Figure~\ref{fig:coll}, a peakon anti-peakon solution $u$ of the Camassa--Holm equation (that is, with $\rho$ identically zero) which enjoys wave breaking. In addition, we have plotted  the corresponding energy function $u^2+u_x^2$. A closer analysis reveals that at the time  $t=t_c$ of wave breaking, all the energy is concentrated at one point, which can be described as (a multiple of) a Dirac delta function. In contrast to that, Figure~\ref{fig:col2} shows that if we choose as initial condition the same peakon anti-peakon function $u_0$ together with $\rho_0(x)=0.5$, then no wave breaking takes place, but at the time $t_c$ where $u(t_c,x)\approx 0$, a considerable part of the energy is transferred from $u^2+u_x^2$  to $\rho^2$, while  the total energy $\int_0^1 (u^2+u_x^2+\rho^2)dx$ remains constant. 

\section{Eulerian and Lagrangian variables}

The two-component Camassa--Holm (2CH) system with $\kappa\in\Real$ and $\eta\in(0,\infty)$ reads
\begin{subequations}\label{eq:chkappa}
  \begin{align}
    u_t-u_{txx}+\kappa u_x+3uu_x-2u_xu_{xx}-uu_{xxx}+\eta\rho\rho_x&=0,\\
    \rho_t+(u\rho)_x&=0.
  \end{align}
\end{subequations}
If $(u,\rho)$ is a solution of \eqref{eq:chkappa} then the pair $(v,\tau)$, given by $v(t,x)=u(t,x-\alpha t)+\alpha$ and $\tau(t,x)=\sqrt{\beta}\rho(t,x)$, is a solutions to the 2CH system with $\kappa$ and $\eta$ replaced by $\kappa-2\alpha$ and $\frac{\eta}{\beta}$, respectively. Thus we can assume without loss of generality that 
$\kappa=0$ and $\eta=1$.  Moreover, we will only consider the Cauchy problem for initial data, and hence also of solutions, of period unity, i.e.,  $u(t,x+1)=u(t,x)$ and 
$\rho(t,x+1)=\rho(t,x)$. All our results carry over with only slight modifications to the case of a general period. 

To any pair $(u_0,\rho_0)$ in $H^1_{\rm per}\times L^2_{\rm per}$ we can introduce the corresponding Lagrangian coordinates $(y(0,\xi),U(0,\xi), \nu(0,\xi), r(0,\xi))$ and describe their time evolution using the weak formulation of the 2CH system. Namely, the characteristics $y(t,\xi)$ are defined as  solutions of 
$$
y_t(t,\xi)=u(t,y(t,\xi))
$$
for a given $y(0,\xi)$ such that $y(0,\xi+1)=y(0,\xi)+1$. The Lagrangian velocity $U(t,\xi)$ defined as 
$$
U(t,\xi)=u(t,y(t,\xi)).
$$ 
The energy derivative reads
$$
\nu(t,\xi)=(u^2+u_x^2+\rho^2)(t,y(t,\xi))y_\xi(t,\xi)
$$
together with the energy $h(t)=\int_0^1 \nu(t,\xi)d\xi$, and, finally, 
$$
r(t,\xi)=\rho(t,y(t,\xi))y_\xi(t,\xi)
$$ 
is the Lagrangian density.

Rewriting the 2CH system as 
\begin{subequations}
  \begin{align}
    u_t+uu_x+P_x&=0,\\
    \rho_t+(u\rho)_x&=0,
  \end{align}
\end{subequations}
where $P=P(t,x)$ implicitly is  given as the solution of $P-P_{xx}=u^2+\frac12 u_x^2+\frac12 \rho^2$, enables us to derive how $(y,U,\nu,r)$ change with respect to time. 
In particular, direct computations yield, after setting $y(t,\xi)=\xi+\zeta(t,\xi)$, that 
\begin{subequations}\label{sys:persys}
  \begin{align}
    \zeta_t&=U,\\
    U_t&=-Q,\\
    \nu_t&=-2QUy_\xi+(3U^2-2P)U_\xi,\\
    r_t&=0,
  \end{align}
\end{subequations}
where 
\begin{align}
  \label{eq:Psimp1}
  P(t,\xi)&=\frac1{2(e-1)}\int_0^1\cosh(y(t,\xi)-y(t,\eta))(U^2y_\xi+\nu)(t,\eta)\,d\eta\\
  &\quad+\frac14\int_0^1\exp\big(-\sgn(\xi-\eta)(y(t,\xi)-y(t,\eta))\big)(U^2y_\xi+\nu)(t,\eta)\,d\eta,\notag
\end{align}
and 
\begin{align}
  \label{eq:Qsimp1}
  Q(t,\xi)&=\frac1{2(e-1)}\int_0^1\sinh(y(t,\xi)-y(t,\eta))(U^2y_\xi+\nu)(t,\eta)\,d\eta \\
  \notag &\quad
  -\frac14\int_0^1\sgn(\xi-\eta)\exp\big(-\sgn(\xi-\eta)(y(t,\xi)-y(t,\eta))\big)(U^2y_\xi+\nu)(t,\eta)\,d\eta.
\end{align}
First, we will consider this system of ordinary differential equations in the Banach space 
$E=W^{1,1}_{\rm per}\times W^{1,1}_{\rm per}\times L^1_{\rm per}\times L^1_{\rm per}$, 
where 
\begin{subequations}
  \begin{align}
    W^{1,1}_{\rm per}&=\{f\in W^{1,1}_{\rm loc}(\Real)\mid f(\xi+1)=f(\xi) \text{ for all }\xi\in\Real\},\\
    L^1_{\rm per}&=\{f\in L^1_{\rm loc}(\Real)\mid f(\xi+1)=f(\xi) \text{ for all }\xi\in\Real\},
  \end{align}
\end{subequations}
and the corresponding norms are given by
\begin{equation*}
  \begin{aligned}
    \norm{f}_{W^{1,1}_{\rm per}}&=\norm{f}_{L^\infty([0,1])}+\norm{f_\xi}_{L^1([0,1])},\text{ and }\norm{f}_{L^1_{\rm per}}=\norm{f}_{L^1([0,1])},\\
    \norm{(y,U,\nu,r)}_{E}&=\norm{y-\id}_{W^{1,1}_{\rm per}}+\norm{U}_{W^{1,1}_{\rm per}}+\norm{\nu}_{L^1_{\rm per}}+\norm{r}_{L^1_{\rm per}}.
  \end{aligned}
\end{equation*}
The existence and uniqueness of short time solutions of \eqref{sys:persys}, will follow from a contraction argument once we can show that the right-hand side of \eqref{sys:persys} is Lipschitz continuous on bounded sets. Note that this is the case if and only if the same holds for $P$ and $Q$. The latter statement  has been proved in \cite[Lemma 2.1]{GHR}, and we state the result here for completeness.

\begin{lemma}
  \label{lem:PQ} 
  For any $X=(y,U,\nu,r)$ in $E$, we define the maps
  $\Q$ and $\PP$ as $\Q(X)=Q$ and $\PP(X)=P$ where
  $P$ and $Q$ are given by \eqref{eq:Psimp1} and
  \eqref{eq:Qsimp1}, respectively. Then, $\PP$ and $\Q$
  are Lipschitz maps on bounded sets from $E$ to
  $\Wperun$. More precisely, we have the following
  bounds. Let 
  \begin{equation}
    \label{eq:defBM}
    B_M=\{X=(y,U,\nu,r)\in E\mid 
    \norm{U}_{\Wperun}+\norm{y_\xi}_{L^1_{per}}+\norm{\nu}_{L^1_{per}}\leq
    M\}.
  \end{equation}
  Then for any $X,\tilde X\in B_M$, we have
  \begin{equation}
    \label{eq:lipbdQ}
    \norm{\Q(X)-\Q(\tilde X)}_{\Wperun}\leq C_M\norm{X-\tilde X}_{E},
  \end{equation}
  and
  \begin{equation}
    \label{eq:lipbdP}
    \norm{\PP(X)-\PP(\tilde X)}_{\Wperun}\leq C_M\norm{X-\tilde X}_{E},
  \end{equation}
  where the constant $C_M$ only depends on the value
  of $M$.
\end{lemma}

To establish the global existence of solutions, we have to impose more conditions on our initial data and solutions in Lagrangian coordinates. 
\begin{definition}
  \label{def:F}
  The set $\F$ is composed of all $(y,U,\nu,r)\in E$
  such that
  \begin{subequations}
    \label{eq:lagcoord}
    \begin{align}
      \label{eq:lagcoord1}
      &\ (y,U)\in
      W^{1,\infty}_{\rm loc}(\Real)\times W^{1,\infty}_{\rm loc}(\Real) ,\ (\nu,r)\in L^\infty(\Real)\times L^\infty(\Real),\\
      \label{eq:lagcoord2}
      &y_\xi\geq0,\ \nu\geq0,\ y_\xi+\nu\geq c\text{
        almost everywhere, for some constant $c>0$},\\
      \label{eq:lagcoord3}
      &y_\xi \nu=y_\xi^2U^2+U_\xi^2+r^2\text{ almost everywhere}.
    \end{align}
  \end{subequations}
\end{definition}

The set $\F$ is preserved with respect to time and plays a special role when proving the global existence of solutions. In particular, for $X(t)\in\F$, we have $h(t)=h(0)$ for all $t\in\Real$ which implies that $\norm{X(t)}_E$ cannot blow up within a finite time interval. Note that the first three equations in \eqref{sys:persys} are independent of $r$ and  coincide with the system considered in \cite{GHR}. Moreover, the last variable $r$ is preserved with respect to time. Hence, by following closely the proofs of \cite[Lemma 2.3, Theorem 2.4]{GHR}, we get the global existence of solutions.
\begin{theorem}
  \label{th:global}
  For any $\bar X=(\bar y,\bar U,\bar \nu,\bar r)\in\F$,
  the system \eqref{sys:persys} admits a unique global
  solution $X(t)=(y(t),U(t),\nu(t),r(t))$ in
  $C^1(\Real_+,E)$ with initial data $\bar X=(\bar
  y,\bar U,\bar\nu,\bar r)$. We have $X(t)\in\F$ for all
  times. Let the mapping
  $S\colon\F\times\Real_+\to\F$ be defined as
  \begin{equation*}
    S_t(X)=X(t).
  \end{equation*}
  Given $M>0$ and $T>0$, we define $B_M$ as before,
  that is,
  \begin{equation}
    \label{eq:defBM2}  
    B_M=\{X=(y,U,\nu,r)\in E\mid 
    \norm{U}_{\Wperun}+\norm{y_\xi}_{L^1_{per}}+\norm{\nu}_{L^1_{per}}\leq
    M\}.
  \end{equation}
  Then there exists a constant $C_M$ which depends
  only on $M$ and $T$ such that, for any two
  elements $X_\alpha$ and $X_\beta$ in $B_M$, we
  have
  \begin{equation}
    \label{eq:stabSt}
    \norm{S_tX_\alpha-S_tX_\beta}_E\leq C_M\norm{X_\alpha-X_\beta}_E
  \end{equation}
  for any $t\in[0,T]$.
\end{theorem}

So far we have proved that there exist global, unique solutions to the 2CH system in Lagrangian coordinates. However, we still have to show that the assumptions are sufficiently general to accommodate rather general initial data in Eulerian coordinates. In particular, we must admit initial data (in Eulerian coordinates) that consists not only of the functions $u_0$ and $\rho_0$ but also of a positive, periodic Radon measure. This is necessary due to the fact that when wave breaking occurs, energy is concentrated at sets of measure zero. More precisely, we define the set of Eulerian coordinates as follows.

\begin{definition}  \label{def:D}
  The set $\D$ of possible initial data consists of all triplets $(u,\rho,\mu)$ such that $u\in H^1_{\rm per}$, $\rho\in L^2_{\rm per}$, and $\mu$ is a positive, periodic Radon measure whose absolute continuous part, $\mu_{\rm ac}$, satisfies 
  \begin{equation}
    \mu_{\rm ac}=(u^2+u_x^2+\rho^2)dx.
  \end{equation}  
\end{definition}

Having identified our set of Eulerian coordinates we can map them to the corresponding set of Lagrangian coordinates, using the mapping $\tilde L$. 

\begin{definition}\label{defL}
  For any $(u,\rho,\mu)$ in $\D$, let
  \begin{equation}\label{def:EtoF}
    \begin{aligned} 
      h& = \mu([0,1)), \\
      y(\xi)& = \sup \{y \mid F_\mu(y)+y<(1+h)\xi\},\\
      \nu(\xi)& = (1+h)-y_\xi(\xi),\\
      U(\xi)& =u\circ y(\xi), \\
      r(\xi)& =\rho\circ y(\xi)y_\xi(\xi)
    \end{aligned}
  \end{equation}
  where 
  \begin{equation}
    F_\mu(x)=\left\{
      \begin{aligned}
        \mu([\,0,x)) \quad & \text{ if } x>0,\\
        0 \quad & \text{ if }x=0,\\
        -\mu([\,x,0)) \quad & \text{ if } x<0.
      \end{aligned}
    \right .
  \end{equation}
  Then $(y,U,\nu,r)\in \F$. We define
  $\tilde L(u,\mu)=(y,U,\nu,r)$. The functions  $P$ and  $Q$   are given by \eqref{eq:Psimp1} and \eqref{eq:Qsimp1}, respectively.
\end{definition}

That this definition is well-posed follows after some slight modifications as in \cite{GHR}. 

However, notice that we have three Eulerian coordinates in contrast to four Lagrangian coordinates,  and hence there can at best be a one-to-one correspondence between triplets in Eulerian coordinates and \textit{equivalence classes} in Lagrangian coordinates. When defining equivalence classes, relabeling functions will play a key role, and we will see why we had to impose \eqref{eq:lagcoord2} in the definition of $\F$. Therefore we will now focus on the set $\Gr$ of relabeling functions.

\begin{definition}\label{def:group}
  \label{def:G} Let $\Gr$ be the set of all
  functions $f$ such that $f$ is invertible,
  \begin{align}
    \label{eq:Gcond1}
    &f\in W_{\rm loc}^{1,\infty}(\Real),\
    f(\xi+1)=f(\xi)+1 \text{ for all }\xi\in\Real\text{, and }\\
    \label{eq:Gcond2}
    &f-\id\text{ and }f^{-1}-\id\text{ both belong to }\Wper.
  \end{align}
\end{definition}
One of the main reasons for the choice of $\Gr$ is that any $f\in \Gr$ satisfies 
\begin{equation*}
  \frac{1}{1+\alpha}\leq f_\xi\leq 1+\alpha,
\end{equation*}
for some constant $\alpha>0$ according to \cite[Lemma 3.2]{HR}. This allows us, following the same lines as in \cite[Definition 3.2, Proposition 3.3]{GHR} to define a group action of $\Gr$ on $\F$.
\begin{definition}
  We define the map $\Phi\colon\Gr\times\F\to\F$
  as follows
  \begin{equation*}
    \left\{
      \begin{aligned}
        \bar y&=y\circ f,\\
        \bar U&=U\circ f,\\
        \bar\nu&=\nu\circ f f_\xi,\\
	\bar r& =r\circ ff_\xi,
      \end{aligned}
    \right.
  \end{equation*}
  where $(\bar y,\bar U,\bar
  \nu,\bar r)=\Phi(f,(y,U,\nu,r))$. We denote $(\bar y,\bar
  U,\bar\nu,\bar r)=(y,U,\nu)\act f$. 
\end{definition}

Using $\Phi$, we can identify a subset of $\F$ which contains one element of each equivalence class.
We introduce $\F_0\subset \F$,
\begin{equation*}
  \F_0=\{ X=(y,U,\nu,r)\in\F\mid y_\xi+\nu=1+h\},
\end{equation*}
where $h=\norm{\nu}_{L^1_{per}}$.
In addition, let $\H\subset \F_0$ be defined as follows
\begin{equation*}
  \H=\{(y,U,\nu,r)\in\F_0\mid \int_0^1 y(\xi)d\xi=0\}.
\end{equation*}
We can then associate to any element $X\in\F$ a unique element $\bar X\in \H$. This means there is a bijection between $\H$ and $\F/\Gr$. Indeed, let 
$\Pi_1\colon\F\to\F_0$ be given by
\begin{equation*}
  \Pi_1(X)=X\act f^{-1},
\end{equation*}
with $f(\xi)=\frac{1}{1+h}(y(\xi)+\int_0^\xi\nu(\eta)d\eta)\in\Gr$ for $X=(y,U,\nu,r)\in \F$ due to \eqref{eq:lagcoord2}
and $\Pi_2\colon\F_0\to\H$ by 
\begin{equation*}
  \Pi_2(x)=X(\xi-a),
\end{equation*}
where $a=\int_0^1 y(\xi)d\xi$ for $X=(y,U,\nu,r)$.
Then
\begin{equation*}
  \Pi=\Pi_2\circ \Pi_1
\end{equation*}
is a projection from $\F$ to $\H$, and, since $\Pi(X)$ is unique, $\F/\Gr$ and $\H$ are in bijection. We can now redefine our mapping from Eulerian to Lagrangian coordinates such that any triplet $(u,\rho,\mu)\in \D$ is mapped to the corresponding element $(y,U,\nu,r)\in \H$ by applying $\tilde L$ followed by $\Pi$. 
\begin{theorem}
  For any $(u,\rho,\mu)\in \D$ let $X=(y,U,\nu,r)\in \H$ be given by $X=L(u,\rho,\mu)=\Pi\circ \tilde L(u,\rho,\mu)$. Then $L\colon\D\to \H$.
\end{theorem}

Furthermore, since $r_t=0$, \cite[Lemma 3.5]{GHR} implies directly that $S_t$ is equivariant, i.e., 
\begin{equation}
  S_t(X\act f)=S_t(X)\act f,
\end{equation}
for $X\in \F$ and $f\in\Gr$. In particular, we can define the semigroup $\bar S_t$ on $\H$ as 
\begin{equation}\label{timeevol}
  \bar S_t=\Pi\circ S_t.
\end{equation}

The final and last step is to go back from Lagrangian to Eulerian coordinates, which is a generalization of \cite[Theorem 3.11]{HR} and the adaptation of \cite[Theorem 4.10]{GHR} to the periodic case. 
\begin{theorem}\label{def:FtoE}
  Let $X\in\F$, then the periodic measure\footnote{The push-forward of a measure $\nu$ by a measurable function $f$ is the measure $f_\#\nu$ defined as $f_\#\nu(B)=\nu(f^{-1}(B))$ for any Borel set $B$.}  $y_\#(r d\xi)$ is absolutely continuous and $(u,\rho,\mu)$ given by
  \begin{equation}
    \label{eq:defL}
    \begin{aligned}
      u(x)=U(\xi) & \text{ for any }   \xi \text{ such that } x=y(\xi),\\
      \mu& =y_\# (\nu d\xi),\\
      \rho(x) dx& =y_\# (r d\xi),
    \end{aligned}
  \end{equation}
  belongs to $\D$. We denote by $M$ the mapping
  from $\F$ to $\D$ which to any $X\in\F$
  associates the element $(u,\rho,\mu)\in\D$ given by
  \eqref{eq:defL}. 
\end{theorem}
Note that the mapping $M$ is independent of the representative in every equivalence class we choose, i.e., 
\begin{equation}
  \label{eq:propM}
  M=M\circ \Pi.
\end{equation}

In order to be able to get back and forth between Eulerian and Lagranian coordinates at any possible time it is left to clarify the relation between $L$ and $M$.

\begin{theorem}\label{thm:LMPI}
  The maps $L:\D\to \H$ and $M:\H\to\D$ are invertible. We have 
  \begin{equation}
    \label{eq:inversma}
    L\circ M=\Pi, \quad \text{ and } \quad M\circ L=\id,
  \end{equation}
  where the mapping $\Pi: \F\to \H$ is a projection which associates to any element $X\in \F$ a unique element $\tilde X\in \H$, which means, in particular, that $\F/\Gr$ and $\H$ are in bijection. 
\end{theorem}

The proof follows the same lines as \cite[Theorem 3.12]{HR}, and we therefore do not present it here. We will see later that the last theorem together with \eqref{timeevol} allows us to define a semigroup of solutions. To obtain a continuous semigroup we have to study the stability of solutions in Lagrangian coordinates, which is the aim of the next section.

\section{Lipschitz metric}

We will now construct a Lipschitz metric in Lagrangian coordinates which will be  invariant under relabeling. It will be quite similar to the one in \cite{GHR} due to the fact that the first three equations in \eqref{sys:persys} are independent of $r$ and coincide with the system considered in \cite{GHR} and because 
$r(t)=r(0)$ for all $t\in\Real$. 

Let $X_\alpha$, $X_\beta\in\F$. We introduce the function $J(X_\alpha, X_\beta)$ by 
\begin{equation}
  \label{eq:defJ}
  J(X_\alpha,X_\beta)=\inf_{f,g\in\Gr}\norm{X_\alpha\act f-X_\beta\act g}_E,
\end{equation}
which is invariant with respect to relabeling.
That means, for any $X_\alpha,X_\beta\in\F$ and
$f,g\in\Gr$, we have
\begin{equation}
  \label{eq:invralJ}
  J(X_\alpha\act f,X_\beta\act g)=J(X_\alpha,X_\beta).
\end{equation}

Note that the mapping $J$ does not define a metric, since it does not satisfy the
triangle inequality, which is the reason why we introduce the
following mapping $d$.

Let  
$d(X_\alpha,X_\beta)$ be defined by 
\begin{equation}
  \label{eq:defdist}
  d(X_\alpha,X_\beta)=\inf \sum_{n=1}^NJ(X_{n-1},X_n), \quad X_\alpha,X_\beta\in\F,
\end{equation}
where the infimum is taken over all finite sequences
$\{X_n\}_{n=0}^N\in\F$ satisfying
$X_0=X_\alpha$ and $X_N=X_\beta$.
In particular, $d$ is relabeling invariant, that means 
for any $X_\alpha,X_\beta\in\F$ and $f,g\in\Gr$,
we have
\begin{equation}
  \label{eq:invralJ_d}
  d(X_\alpha\act f,X_\beta\act g)=d(X_\alpha,X_\beta).
\end{equation}

In order to prove that $d$ is a Lipschitz metric on bounded sets, we have to choose one element in each equivalence class, and we will apply \eqref{eq:stabSt}. One problem we are facing in that context is that the constant on the right-hand side of \eqref{eq:stabSt} depends on the set $B_M$ we choose, but $B_M$ is not preserved by the time evolution while it is invariant with respect to relabeling. Hence we will try to find a suitable set, which is invariant with respect to time and relabeling and is in some sense equivalent to $B_M$. 
To that end we define the subsets of bounded energy $\F^M$
of $\F_0$ by
\begin{equation*}
  \F^M=\{X=(y,U,\nu,r)\in \F\mid  h=\norm{\nu}_{L^1_{per}}\leq M\}
\end{equation*}
and let $\H^M=\H\cap\F^M$.
The important property of the set $\F^M$ is that it
is preserved both by the flow and relabeling. In particular, 
we have that
\begin{equation}
  \label{eq:incFMBM}
  B_{M}\cap\H\subset \H^M \subset B_{\bar M}\cap\H
\end{equation}
for $\bar M=6(1+M)$ and hence the sets $B_M\cap\H$
and $\H^M$ are in this sense equivalent.

\begin{definition}
  Let $d_M$ be the metric on $\H^M$ which is
  defined, for any $X_\alpha,X_\beta\in\H^M$, as
  \begin{equation}
    \label{eq:defdM}
    d_M(X_\alpha,X_\beta)=\inf \sum_{n=1}^NJ(X_{n-1},X_n)
  \end{equation}
  where the infimum is taken over all finite
  sequences $\{X_n\}_{n=0}^N\in\H^M$ which
  satisfy $X_0=X_\alpha$ and $X_N=X_\beta$.
\end{definition}

By definition $d_M$ is relabeling invariant and the triangle inequality is satisfied. In this way we obtain a metric which in addition can be compared with other norms on $\H^M$ (cf.~\cite[Lemma 4.3]{GHR}).
\begin{lemma}
  The mapping $d_M\colon\H^M\times \H^M\to \Real_+$ is a metric on $\H^M$. Moreover,  
  given $X_\alpha$, $X_\beta\in \H^M$, define $R_\alpha=\int_0^1 r_\alpha(\eta)d\eta$ and $R_\beta=\int_0^1 r_\beta(\eta)d\eta$. Then we have 
  \begin{equation}
    \label{eq:LinfbdJ}
    \norm{y_\alpha-y_\beta}_{L^\infty}+\norm{U_\alpha-U_\beta}_{L^\infty}+\abs{h_\alpha-h_\beta}+\abs{R_\alpha-R_\beta}\leq C_M d_M(X_\alpha,X_\beta)
  \end{equation}
  and 
  \begin{equation}
    \label{eq:dequiv}
    d(X_\alpha,X_\beta)\leq\norm{X_\alpha-X_\beta}_E,
  \end{equation}
  where $C_M$ denotes some fixed constant which depends only on $M$.
\end{lemma}

To show that we not only obtained a relabeling invariant metric but in fact a Lipschitz metric, we combine all results we obtained so far as in \cite[Theorem 4.6]{GHR}. This yields the following Lipschitz stability theorem for $\bar S_t$.

\begin{theorem}
  \label{th:stab} Given $T>0$ and $M>0$, there
  exists a constant $C_M$ which depends only on $M$
  and $T$ such that, for any
  $X_\alpha,X_\beta\in\H^M$ and $t\in[0,T]$, we
  have
  \begin{equation}
    \label{eq:stab}
    d_M(\bar S_tX_\alpha,\bar S_tX_\beta)\leq C_Md_M(X_\alpha,X_\beta).
  \end{equation}
\end{theorem}

\section{Global weak solutions}

It is left to check that we obtain a global weak
solution of the 2CH system by solving
\eqref{sys:persys} and using the maps between
Eulerian and Lagrangian coordinates.  In the case
of conservative solutions we have that for any
triplet $(u(t,x),\rho(t,x),\mu(t,x))$ in Eulerian
coordinates, the function $P(t,x)$ is given by
\begin{equation}\label{Peuler}
  \begin{aligned}
    P(t,x)&=\frac{1}{2(e-1)}\int_0^1 \cosh(x-z)u^2(t,z)dz+\frac{1}{2(e-1)}\int_0^1 \cosh(x-z)d\mu(t,z)\\
    &\quad +\frac14 \int_0^1e^{-\vert x-z\vert}u^2(t,z)dz+\frac14 \int_0^1 e^{-\vert x-z\vert}d\mu(t,z).
  \end{aligned}
\end{equation}
Applying the mapping $L$ maps $P(t,x)$ to
$P(t,\xi)$ given by \eqref{eq:Psimp1} and
$P_x(t,x)$ to $Q(t,\xi)$ given by
\eqref{eq:Qsimp1}.  Since the set of times
where wave breaking occurs has measure zero,
$\mu=\mu_{ac}=(u^2+u_x^2+\rho^2)dx$ for almost
all times, $P(t,x)$ defined by \eqref{Peuler}
coincides for almost all times and all
$x\in\Real$ with the solution of
$P-P_{xx}=u^2+\frac12 u_x^2+\frac12 \rho^2$.

\begin{definition}\label{weaksol}
  Let $u\colon\Real_+\times\Real \rightarrow
  \Real$ and
  $\rho\colon\Real_+\times\Real\rightarrow\Real$.
  Assume that $u$ and $\rho$
  satisfy \\
  (i) $u\in L^\infty([0,\infty), H^1_{\rm per})$, $\rho\in L^\infty([0,\infty), L^2_{\rm per})$, \\
  (ii) the equations
  \begin{multline}\label{weak1}
    \iint_{\Real_+\times [0,1]}\Big(-u(t,x)\phi_t(t,x)+(u(t,x)u_x(t,x)+P_x(t,x))\phi(t,x)\Big)dxdt \\
    = \int_{[0,1]} u(0,x)\phi(0,x)dx,
  \end{multline}
  \begin{equation}\label{weak2}
    \iint_{\Real_+\times [0,1]}\Big((P(t,x)-u^2(t,x)-\frac{1}{2} u_x^2(t,x))\phi(t,x)+P_x(t,x)\phi_x(t,x)\Big)dxdt=0,
  \end{equation}
  and \begin{equation}\label{weak3}
    \iint_{\Real_+\times [0,1]}\Big(-\rho(t,x)\phi_t(t,x)-u(t,x)\rho(t,x)\phi_x(t,x)\Big) dxdt=\int_{[0,1]} \rho(0,x)\phi(0,x)dx,
  \end{equation}
  hold for all spatial periodic functions $\phi\in
  C_0^\infty ([0,\infty),\Real)$.  Then we say
  that $(u,\rho)$ is a global weak solution of the
  two-component Camassa--Holm system.
  
  If
  $(u,\rho)$ in addition satisfies
  \begin{equation*}
    (u^2+u_x^2+\rho^2)_t+(u(u^2+u_x^2+\rho^2))_x-(u^3-2Pu)_x=0
  \end{equation*}
  in the sense that
  \begin{align}
    \iint_{\Real_+\times [0,1]}&\Big[
    (u^2(t,x)+u_x^2(t,x)+\rho^2(t,x))\phi_t(t,x)\notag\\
    &+(u(t,x)(u^2(t,x)+u_x^2(t,x)+\rho^2(t,x)))\phi_x(t,x)\label{eq:weak4}\\
    &\qquad\qquad-(u^3(t,x)-2P(t,x)u(t,x))\phi_x(t,x)\Big]dxdt \notag
    =0,
  \end{align}
  for any spatial periodic function $\phi\in
  C_0^{\infty}((0,\infty)\times\Real)$, we say
  that $(u,\rho)$ is a weak global conservative
  solution of the two-component Camassa--Holm
  system.
\end{definition}

Introduce the mapping $T_t$ from $\D$ to $\D$ by 
\begin{equation}
  T_t=M\bar S_t L. 
\end{equation}
Then one can check that for any $(u_0,\rho_0,\mu_0)\in \D$ such that $\mu_0$ is purely absolutely continuous, the pair $(u(t,x),\rho(t,x))$ given by $(u,\rho,\mu)(t)=T_t(u_0,\rho_0,\mu_0)$ satisfies \eqref{weak1}--\eqref{weak3}. 
\begin{theorem}\label{th:main2}
  Given any initial condition $(u_0,\rho_0)\in H^1_{\rm per}\times L^2_{\rm per}$, we define $\mu_0=(u_0^2+u_{0,x}^2+\rho_0^2)dx$, and we denote $(u,\rho,\mu)(t)=T_t(u_0,\rho_0, (u_0^2+u_{0,x}^2+\rho_0^2)dx)$. Then $(u,\rho)$ is a periodic and global weak solution of the 2CH system and $\mu$ satisfies weakly $$\mu_t+(u\mu)_x=(u^3-2Pu)_x.$$ Moreover, $\mu(t)$ consists of an absolutely continuous and a singular part, that means 
  \begin{equation*}
    \int_{[0,1]} d\mu(t,x)=\int_{[0,1]} (u^2+u_{x}^2+\rho^2)(t,x) \, dx+\int_{[0,1]} d\mu_{sing}(t,x).
  \end{equation*}
  In particular, $supp(\mu_{sing}(t,.))$ coincides with the set of points where wave breaking occurs at time $t$.
\end{theorem}

Note that since we are looking for global weak solutions for initial data in $H^1_{\rm per}\times L^2_{\rm per}$ it is no restriction to assume that $\mu$ is purely absolutely continuous initially  while we in general will not have that $\mu$ remains purely absolutely continuous at any later time. In particular, if $\mu$ is not absolutely continuous at a particular time, we know how much and where the energy has concentrated, and this energy must be given back to the solution in order to obtain conservative solutions. Therefore the measure plays an important role. \\
It is also possible to define global weak solutions for initial data where the measure $\mu_0$ is not purely absolutely continuous by defining $P(0,x)$ using \eqref{Peuler}, which is then mapped to \eqref{eq:Psimp1} by applying \eqref{def:EtoF} directly. Moreover, $P(0,\xi)$ can be mapped back to $P(0,x)$ via $M$. In addition, this point of view allows us to jump between Eulerian and Lagrangian coordinates at any time. \\
Moreover, one can show that the sets $\F^M$ and  $\H^M$ in Lagrangian coordinates correspond to the set $\D^M$ in Eulerian coordinates. Given $M>0$, we define 
\begin{equation}
  \D^M=\{(u,\rho,\mu)\in\D\mid \mu([0,1))\leq M\}.
\end{equation}
Thus it is natural to define a Lipschitz metric on the sets of bounded energy in Eulerian coordinates as follows,
\begin{equation}
  d_{\D^M}((u,\rho,\mu), (\tilde u,\tilde \rho,\tilde\mu))=d_M(L(u,\rho,\mu), L(\tilde u,\tilde\rho,\tilde\mu)).
\end{equation}
In particular, we have the following result. 

\begin{theorem}\label{th:main}
  The semigroup $(T_t,d_\D)$, which corresponds to solutions of the 2CH system, is a continuous
  semigroup on $\D$ with respect to the metric
  $d_\D$. The semigroup is Lipschitz continuous on
  sets of bounded energy, that is, given $M>0$ and
  a time interval $[0,T]$, there exists a constant
  $C$ which only depends on $M$ and $T$ such that,
  for any $(u,\rho,\mu)$ and $(\tilde u,\tilde \rho,\tilde\mu)$ in
  $\D^M$, we have
  \begin{equation*}
    d_{\D^M}(T_t(u,\rho,\mu),T_t(\tilde{u},\tilde\rho,\tilde{\mu}))\leq Cd_{\D^M}((u,\rho,\mu),(\tilde{u},\tilde\rho,\tilde{\mu}))
  \end{equation*}
  for all $t\in[0,T]$.
\end{theorem}

Last, but not least, we want to investigate the regularity of solutions and the connection of the topology in $\D$ with other topologies. Due to the global interaction term given for almost all times by 
\begin{equation}\label{eq:nonlocal}
  \begin{aligned}
    P(t,x)& =\frac1{2(e-1)} \int_0^1 \cosh(x-z)(2u^2+u_x^2+\rho^2)(t,z)dz\\ 
    &\quad + \frac14 \int_0^1 e^{-\vert x-z\vert}(2u^2+u_x^2+\rho^2)(t,z)dz, 
  \end{aligned}
\end{equation}
the 2CH system has an infinite speed of propagation \cite{henry:09}. However, the system remains
essentially hyperbolic in nature, and we prove that
singularities travel with finite speed.  In \cite[Theorem 6.1]{GHR2} we showed that the local regularity of a solution depends on the regularity of the initial data and that  $\rho_0(x)^2$ can be bounded from below by a strictly positive constant. Since this result is a local result, it carries over to the periodic case and we state it here for the sake of completeness.

\begin{theorem}\label{th:presreg}
  We consider initial data $(u_0,\rho_0,\mu_0)\in\D$. Furthermore, we assume that there exists an interval $(x_0,x_1)$ such that $(u_0,\rho_0,\mu_0)$ is $p$-regular, with $p\geq 1$, in the sense that
  \begin{equation*}
    u_0\in W^{p,\infty}(x_0,x_1), \quad \rho_0\in W^{p-1,\infty}(x_0,x_1), \quad \text{and } \mu_0=\mu_{0,ac} \text{ on } (x_0,x_1), 
  \end{equation*}
  and that 
  \begin{equation*}
    \rho_0(x)^2\geq c>0
  \end{equation*}
  for $x\in(x_0,x_1)$. Then for any $t\in\Real_+$, $(u,\rho,\mu)(t,\dott)$ is $p$-regular on the interval $(y(t,\xi_0),y(t,\xi_1))$, where $\xi_0$ and $\xi_1$ satisfy $y(0,\xi_0)=x_0$ and $y(0,\xi_1)=x_1$ and are defined as 
  \begin{equation*}
    \text{$\xi_0=\sup\{\xi\in\Real\ |\ y(0,\xi)\leq x_0\}$  and 
      $\xi_1=\inf\{\xi\in\Real\ |\ y(0,\xi)\geq x_1\}$.}
  \end{equation*}
\end{theorem}

In other words, we see that the regularity is preserved between characteristics. As an immediate consequence we obtain the following result.

\begin{theorem}\label{cor:presreg3}
  If the initial data $(u_0,\rho_0,\mu_0)\in\D$
  satisfies $u_0,\rho_0\in C^{\infty}(\Real)$,
  $\mu_0$ is absolutely continuous and
  $\rho_0^2(x)\geq d>0$ for all $x\in\Real$, then
  $u,\rho\in C^\infty(\Real\times\Real)$
  is the unique classical solution to
  \eqref{eq:chkappa} with 
  $\kappa=0$ and $\eta=1$.
\end{theorem}  

In particular this result implies that if $\rho_0^2(x)\geq c$ for some positive constant $c>0$, no wave breaking occurs. Hence if we can compare the topology on $\D$ with standard topologies we have a chance to approximate conservative solutions of the CH equation which enjoy wave breaking by global smooth solutions of the 2CH system. Indeed,
the mapping 
\begin{equation}
  (u,\rho)\mapsto (u,\rho, (u^2+u_x^2+\rho^2)dx),
\end{equation}
is continuous from $H^1_{\rm per}\times L^2_{\rm per}$ to $\D$. This means, given a sequence $(u_n,\rho_n)\in H^1_{\rm per}\times L^2_{\rm per}$ converging to $(u,\rho)\in H^1_{\rm per}\times L^2_{\rm per}$, then $(u_n, \rho_n, (u_n^2+u_{n,x}^2+\rho_n^2)dx)$ converges to $(u,\rho, (u^2+u_x^2+\rho^2)dx)$ in $\D$. \\
Conversely if $(u_n,\rho_n,\mu_n)$ is a sequence in $\D$ which converges to $(u,\rho,\mu)\in\D$, then 
\begin{equation}
  u_n\rightarrow u \text{ in } L^\infty_{\rm per},\quad \rho_n\overset{\ast}{\rightharpoonup}\rho, \text{ and } \mu_n \overset{\ast}{\rightharpoonup}\mu.  
\end{equation}

Putting now everything together we have the following result.
\begin{theorem}
  \label{th:approxCH}
  Let $u_0\in H^1_{\rm per}$. We consider the
  approximating sequence of initial data
  $(u_0^n,\rho_0^n,\mu_{0}^n)\in\D$ given by
  $u_0^n\in C^\infty(\Real)$ with
  $\lim_{n\to\infty}u_0^n=u_0$ in
  $H^1_{\rm per}$, $\rho_0^n\in C^\infty(\Real)$
  with $\lim_{n\to\infty}\rho_0^n=0$ in
  $L^2_{\rm per}$, $(\rho_0^n)^2\geq d_n$ for
  some constant $d_n>0$ and for all $n$ and
  $\mu_{0}^n=((u_{0,x}^{n})^2+(\rho_0^n)^2)\,dx$. We
  denote by $(u^n,\rho^n)$ the unique classical
  solution to \eqref{eq:chkappa}, with   $\kappa=0$ and $\eta=1$, in
  $C^\infty(\Real_+\times\Real)\times
  C^\infty(\Real_+\times\Real)$ with $(u,\rho)|_{t=0}=(u_0,\rho_0)$.  Then for every
  $t\in\Real_+$, the sequence $u^n(t,\dott)$
  converges to $u(t,\dott)$ in
  $L^{\infty}(\Real)$, where $u$ is the
  conservative solution of the Camassa--Holm
  equation 
  \begin{equation}
    u_t-u_{txx}+3uu_x-2u_xu_{xx}-uu_{xxx}=0,
  \end{equation}
  with initial data $u_0\in
  H^1_{\rm per}$.
\end{theorem}


\begin{thebibliography}{XXX}


\bibitem{BreCons:07}
  A. Bressan and A. Constantin.
  \newblock Global conservative solutions of the {C}amassa--{H}olm equation.
  \newblock {\em Arch. Ration. Mech. Anal.}, 183(2):215--239, 2007.

\bibitem{BreCons:09}
  A. Bressan and A. Constantin.
  \newblock Global dissipative solutions of the Camassa--Holm equation. 
  \newblock {\em Analysis and Applications}, 5:1--27, 2007.

\bibitem{BHR}
  A. Bressan, H. Holden, and X. Raynaud.
  \newblock Lipschitz metric for the Hunter--Saxton equation.
  \newblock {\em J. Math. Pures Appl.}, 94:68--92, 2010.
  
\bibitem{CH:93}
  R.~Camassa and D.~D. Holm.
  \newblock An integrable shallow water equation with peaked solitons.
  \newblock {\em Phys. Rev. Lett.}, 71(11):1661--1664, 1993.

\bibitem{ChenLiu2010}
  R. M. Chen and Y. Liu.
  \newblock Wave breaking and global existence for a generalized two-component Camassa--Holm system.
  \newblock {\em Inter. Math Research Notices}, Article ID rnq118, 36 pages, 2010.

\bibitem{MR2474608}
  A. Constantin and R.~I. Ivanov.
  \newblock On an integrable two-component {C}amassa--{H}olm shallow water
  system.
  \newblock {\em Phys. Lett. A}, 372(48):7129--7132, 2008.


\bibitem{eschlechyin:07} 
  J. Escher, O. Lechtenfeld, and Z. Yin. 
  \newblock Well-posedness and blow-up phenomena for the 2-component
  {C}amassa--{H}olm equation.
  \newblock {\em Discrete Contin. Dyn. Syst.}, 19(3):493--513, 2007.

\bibitem{FuQu:09}
  Y. Fu and C. Qu.
  \newblock Well posedness and blow-up solution for a new coupled 
  {C}amassa--{H}olm equations with peakons.
  \newblock {\em J. Math. Phys.}, 50:012906, 2009.

\bibitem{GHRb:10} K. Grunert, H. Holden, and  X. Raynaud. 
  \newblock Lipschitz metric for the  {C}amassa--{H}olm equation on the  line. 
  \newblock  {\em Discrete Contin. Dyn. Syst.}, (to appear).
  

\bibitem{GHR:12} K. Grunert, H.  Holden, and X.  Raynaud.  
  \newblock Global conservative  solutions of the {C}amassa--{H}olm equation for
  initial data with nonvanishing asymptotics. 
  \newblock {\em  Discrete Contin. Dyn. Syst.}, 32:4209--4227, 2012.
  

\bibitem{GHR}
  \newblock K.~Grunert, H.~Holden, and X.~Raynaud.
  \newblock Lipschitz metric for the periodic Camassa--Holm equation. 
  \newblock {\em J. Differential Equations} 250: 1460--1492, 2011. 

\bibitem{GHR2}
  \newblock K.~Grunert, H.~Holden, and X.~Raynaud.
  \newblock Global solutions for the two-component Camassa--Holm system.
  \newblock {\em Comm. Partial Differential Equations}, 37:2245--2271, 2012.



\bibitem{GuanYin2010} 
  C. Guan and Z. Yin.
  \newblock Global weak solutions for a modified two-component Camassa--Holm equation.
  \newblock {\em Ann. I. H. Poincar{\'e} -- AN} 28:623--641, 2011.

\bibitem{GuanYin2010a} 
  C. Guan and Z. Yin.
  \newblock Global existence and blow-up phenomena for an integrable  two-component Camassa--Holm shallow water system.
  \newblock {\em J. Differential Equations},  248:2003--2014, 2010.

\bibitem{GuiLiu2011}  
  G. Gui and Y. Liu.
  \newblock On the Cauchy problem for the two-component Camassa--Holm system.
  \newblock {\em Math Z},  268:45--66, 2011.

\bibitem{GuiLiu2010}  
  G. Gui and Y. Liu.
  \newblock On the global existence and wave breaking criteria for the two-component Camassa--Holm system.
  \newblock {\em J. Func. Anal.},  258:4251--4278, 2010.

\bibitem{GuoZhou2010}
  Z. Guo and Y. Zhou.
  \newblock On solutions to a two-component generalized Camassa--Holm equation.
  \newblock {\em Studies Appl. Math.}, 124:307--322, 2010.

\bibitem{henry:09}
  D. Henry.
  \newblock Infinite propagation speed for a two component {C}amassa--{H}olm
  equation.
  \newblock {\em Discrete Contin. Dyn. Syst. Ser. B}, 12(3):597--606, 2009.

\bibitem{HolRay:07}
  H. Holden and X. Raynaud.
  \newblock Global conservative solutions of the {C}amassa--{H}olm equation---a
  {L}agrangian point of view.
  \newblock {\em Comm. Partial Differential Equations}, 32(10-12):1511--1549,
  2007.

\bibitem{HolRay:09}
  H. Holden and X. Raynaud.
  \newblock Dissipative solutions for the Camassa--Holm equation. 
  \newblock {\em Discrete Contin. Dyn. Syst.} 24:1047--1112, 2009.


\bibitem{HR}
  \newblock H.~Holden and X.~Raynaud.
  \newblock Periodic conservative solutions of the Camassa--Holm equation.
  \newblock {\em Ann. Inst. Fourier (Grenoble)} 58:945--988, 2008.


\bibitem{Kuzmin} 
  P. A.  Kuz'min.
  \newblock Two-component generalizations of the Camassa--Holm equation.
  \newblock {\em Math. Notes},  81:130--134, 2007.

\bibitem{OlverRosenau}
  P.~J. Olver and P. Rosenau.
  \newblock Tri-hamiltonian duality between solitons and solitary-wave solutions
  having compact support.
  \newblock {\em Phys. Rev. B}, 53(2):1900--1906, 1996.

\bibitem{WangHuangChen}  
  Y. Wang,  J. Huang, and L. Chen.
  \newblock Global conservative solutions of the two-component Camassa--Holm shallow water system.
  \newblock {\em Int.  J.  Nonlin. Science},  9:379--384, 2009.


\end{thebibliography}
\end{document}